\documentclass[12pt]{article}

\begin{document}

\newtheorem{thm}{Theorem}[section]
\newtheorem{thmA}[thm]{Theorem A}
\newtheorem{cor}[thm]{Corollary}
\newtheorem{prop}[thm]{Proposition}
\newtheorem{define}[thm]{Definition}
\newtheorem{rem}[thm]{Remark}
\newtheorem{example}[thm]{Example}
\newtheorem{lemma}[thm]{Lemma}
\def\theequation{\thesection.\arabic{equation}}

\title{Complex-valued Burgers and KdV-Burgers equations}

\author{Netra Khanal\thanks{Department of Mathematics, Oklahoma State
University, 401 Mathematical Sciences, Stillwater, OK 74078, USA.
Email: nkhanal@math.okstate.edu.} \and Jiahong Wu\thanks{Department
of Mathematics, Oklahoma State University, 401 Mathematical
Sciences, Stillwater, OK 74078, USA. Email:
jiahong@math.okstate.edu.} \and Juan-Ming Yuan \thanks{Department of
Applied Mathematics, Providence University, Shalu, Taichung 433,
Taiwan. Email: jmyuan@pu.edu.tw} \and Bing-Yu
Zhang\thanks{Department of Mathematical Sciences, University of
Cincinnati,   Cincinnati, OH 45221-0025, USA. Email:
zhangb@email.uc.edu}}
\date{}
\maketitle

\begin{abstract}
Spatially periodic complex-valued solutions of the Burgers and
KdV-Burgers equations are studied in this paper. It is shown that
for any sufficiently large time $T$, there exists an explicit
initial data such that its corresponding solution of the Burgers
equation blows up at $T$. In addition, the global convergence and
 regularity of   series solutions is established for initial data
satisfying mild conditions.
\end{abstract}

\begin{center}
{\bf Key Words:}

\vspace{.1in} Complex Burgers equation, complex KdV-Burgers
equation, finite-time singularity, global regularity

\vspace{.2in}
{\bf Mathematics Subject Classification:}

\vspace{.1in}
35Q53, 35B35, 35B65, 76D03
\end{center}

\newpage
\section{Introduction}
This work addresses the global regularity issue on
solutions of the complex Burgers and KdV-Burgers equations
\begin{equation}\label{kdv-burgers}
u_t -6 uu_x + \alpha u_{xxx} -\nu u_{xx} =0,
\end{equation}
where $\nu \ge 0$ and $\alpha \ge 0$ are parameters and $u=u(x,t)$
is a complex-valued function. Attention will be focused on the
spatially periodic solutions, namely  $x\in {\mathbf T}= {\mathbf R}/(2\pi),$
the one-torus and we supplement (\ref{kdv-burgers}) with a given
initial data
\begin{equation}\label{ini}
u(x,0) =u_0(x), \quad x\in {\mathbf T}.
\end{equation}

Our first major result is for the complex Burgers equation
((\ref{kdv-burgers}) with $\alpha=0$) and it asserts that for any
sufficiently large time $T$, there exists an explicit smooth initial
data $u_0$ such that its corresponding solution blows up at $t=T$
(Theorem \ref{bup}). This result was partially motivated by a recent
paper of Pol\'{a}\v{c}ik and \v{S}ver\'{a}k \cite{PS}, in which the
complex-valued Burgers equation on the whole line was shown to
develop finite-time singularities for compactly supported smooth
data. Their proof takes advantage of the explicit solution formula
obtained via the Hopf-Cole transform. By contrast, the finite-time
singular solutions constructed in this paper assume the form
\begin{equation}\label{sss}
u(x,t) = \sum_{k=1}^\infty a_k(t)\, e^{ikx}
\end{equation}
and correspond to the initial data $u_0(x)=a\,e^{ix}$. We emphasize
that solutions of the form (\ref{sss}) are locally well-posed in the
usual Sobolev space $H^s:=H^s({\mathbf T})$ with a suitable index
$s$ (see Theorem \ref{local1} for more details). For any $T\ge T_0$
(a fixed number depending on $\nu$ only), we obtain a lower bound
for $|a_k(T)|$ through a careful observation of the pattern that
$a_k(t)$'s exhibit and the finite time singularity of (\ref{sss}) in
$L^2$ then follows if we take $a$ in $u_0$ to be sufficiently large.
This result reveals a fundamental difference between the real-valued
solutions of the Burgers equation and their complex counterparts.
The diffusion in the case of complex-valued solutions no longer
dissipates the $L^2$-norm, which can blow up in a finite time.
However, if we know the $L^2$-norm of a complex-valued solution is
bounded, then there would be no finite-time singularity (Theorem
\ref{nu}).

\vspace{.1in} We also explore the conditions under which solutions
of (\ref{kdv-burgers}) are global in time. A simple example of the
global solutions of (\ref{kdv-burgers}) corresponds to the initial
data $u_0(x)=a_0 e^{ix}$ with $|a_0|<1$ provided $\nu$ and $\alpha$
satisfy a suitable condition, say $\nu^2 + 4 \alpha^2 \ge 9$ (see
Theorem \ref{sing}). For general initial data of the form
$$
u_0(x)  =\sum_{k=1}^\infty a_{0k}\,e^{ikx}
$$
with $|a_{0k}| <1$, (\ref{kdv-burgers}) possesses a unique local
solution (\ref{sss}) with $a_k(t)$ given by a finite sum of terms
that can be made explicit through an inductive relation. To show the
convergence of (\ref{sss}) for large time, it is necessary to
estimate $|a_k(t)|$ and our approach is to count the total number of
terms that it contains. This counting problem is closely related to
the number of nonnegative integer solutions to the equation
$$
j_1 + 2j_2 +3j_3 +\cdots +kj_k =k
$$
for a fixed integer $k>0$. Using a result by Hardy and Ramanujan \cite{HR},
we are able to establish the global regularity of (\ref{sss}) under a mild assumption
(see Theorem \ref{series2}). In addition, $\|u(\cdot,t)\|_{H^s}$ for any $s\ge 0$
decays exponentially in $t$ for large $t$.


\vspace{.1in} We remark that the study of complex-valued Burgers and
KdV-Burgers equations can be justified both physically and
mathematically. Physically these complex equations do arise in the
modeling of several physical phenomena
(\cite{Ke},\cite{Le},\cite{LeSa}). Mathematically these equations
exhibit some remarkable features and admit solutions with much
richer structures than those of their real-valued ones. In fact,
these equations and other complex-valued partial differential
equations have attracted quite some attention recently. A lot of
efforts have been devoted to the important issue of whether or not
their solutions can blow up in a finite time. In \cite{Bir} Birnir
considered the complex KdV equation and constructed a family of
singular solutions represented by the Weierstrass function. Very
recently  Y. Li \cite{Li} obtained simple explicit formulas for
finite-time blowup solutions of the complex KdV equation through
Darboux transform. In \cite{BW} Bona and Weissler addressed the
blowup issue for a family of complex-valued nonlinear dispersive
equations. The papers of Yuan and Wu
(\cite{WY1},\cite{WY2},\cite{YW1}) treated the complex KdV and
KdV-Burgers equations as systems of two nonlinearly coupled
equations and clarified how the potential singularities of the real
part are related to those of the imaginary part. In addition,
extensive numerical experiments were performed to reveal the blowup
structures. Another important example that shows significant
differences between the real-valued and complex-valued solutions is
the Navier-Stokes equations. It remains open whether or not
classical solutions of the 3D incompressible Navier-Stokes equations
can develop finite-time singularities. However, Li and Sinai
\cite{LS} recently showed that the complex solutions of the 3D
Navier-Stokes equations corresponding to large parameter family of
initial data blow up in finite time. Their work motivated the study
of Pol\'{a}\v{c}ik and \v{S}ver\'{a}k on the complex-valued
solutions of the Burgers equation, as we mentioned earlier.

\vspace{.1in} The rest of this paper is divided into three sections.
The second section focuses on the complex Burgers equation and
presents Theorems \ref{bup}, \ref{local1} and \ref{nu}. The third
section details the global regularity results concerning the complex
KdV-Burgers equations.

\vspace{.1in}
\section{Blowup for the complex Burgers equation}
\setcounter{equation}{0}
\label{sec:blow}

This section presents three major results. The first one is a blowup result for the complex
Burgers equation in a periodic domain ${\mathbf T}$$=[0,2\pi]$,
namely
\begin{equation}\label{ivpbur}
\left\{
\begin{array}{l}
 u_t -6 uu_x -\nu u_{xx} =0,\quad x\in {\mathbf T},\,\,t>0, \\
u(x,0)=u_0(x), \quad x\in {\mathbf T}
\end{array}
\right.
\end{equation}
It states that for any sufficiently large $T>0$, there exists an initial data $u_0$ such that
its corresponding solution $u$ blows up at $t=T$. This solution can be represented by
\begin{equation}\label{form}
u(x,t) = \sum_{k=1}^\infty a_k(t)\, e^{ikx}
\end{equation}
and the blowup is in the $L^2$ sense.

\vspace{.1in}
For the sake of completeness of our theory on (\ref{ivpbur}), we also present a local existence and uniqueness
result on solutions of the form (\ref{form}) to the complex-valued KdV-Burgers type equation
\begin{equation}\label{kdvb}
u_t - 6 uu_x + \nu (-\Delta)^\gamma u + \alpha\, u_{xxx} =0,
\end{equation}
which reduces to the complex Burgers equation when $\gamma=1$ and $\alpha=0$. The fractal
Laplacian $(-\Delta)^\gamma$ is defined through Fourier transform,
$$
\widehat{(-\Delta)^\gamma u} (\xi) = |\xi|^{2\gamma}\, \widehat{u}(\xi).
$$
The third result asserts that if the $L^2$-norm of a solution of (\ref{kdvb}) is bounded on $[0,T]$,
then all higher derivatives are bounded and no singularity is possible on $[0,T]$.

\vspace{.1in}
We divide the rest of this section into two subsections with the first devoted to the blowup result and
the second to the local existence uniqueness.

\vspace{.2in}
\noindent{\bf \ref{sec:blow}.1. Finite-time blowup}

\vspace{.1in}
\begin{thm}\label{bup}
For every sufficiently large $T>0$, there exists an initial data $u_0$ of the form
\begin{equation}\label{u00}
u_0(x) = a \,e^{ix}
\end{equation}
such that the corresponding solution $u$ of (\ref{ivpbur}) blows up at $t=T$ in the $L^2$-norm, namely
\begin{equation}\label{l2b}
\|u(\cdot,T)\|_{L^2({\mathbf T})} = \infty.
\end{equation}
\end{thm}

\vspace{.1in}
For any $s\in {\mathbf R}$, the homogeneous Sobolev space $\mathring{H}^s({\mathbf T})$ and the inhomogeneous
Sobolev space $H^s({\mathbf T})$ are defined in the standard fashion. In particular, a function of the form
$$
u(x,t) = \sum_{k=1}^\infty \,a_k \,e^{ikx}
$$
is in $\mathring{H}^s({\mathbf T})$  if
$$
\|u\|^2_{\mathring{H}^s({\mathbf T})}  \equiv \sum_{k=1}^\infty k^{2s}\, |a_k|^2 < \infty,
$$
and in $H^s({\mathbf T})$ if
$$
\|u\|^2_{H^s({\mathbf T})}  \equiv \sum_{k=1}^\infty (1+k^2)^s\, |a_k|^2 < \infty.
$$
Clearly, $L^2({\mathbf T})$ can be identified with $H^0({\mathbf T})$.

\vspace{.1in}
For $u_0$ given by (\ref{u00}), the local existence and uniqueness result of the next subsection asserts
that the corresponding solution $u$ can be written as
$$
u(x,t) = \sum_{k=1}^\infty a_k(t)\, e^{ikx}
$$
before it blows up. The idea is to choose large $a$ such that
$$
\|u(\cdot,T)\|_{L^2}^2 = \sum_{k=1}^\infty |a_k(T)|^2 = \infty.
$$
We attempt to find an explicit representation for $a_k(t)$. It is easy to verify the following iterative formula
\begin{equation}\label{iter}
a_1(t)=a\,e^{-\nu t},\quad a_k(t) = 3ik\, e^{-\nu k^2\,t} \int_0^t e^{\nu k^2 \,\tau}
\sum_{k_1+k_2=k} a_{k_1}(\tau)\, a_{k_2}(\tau) \, d\tau, \quad k=2,3,\cdots.
\end{equation}
To see the pattern in $a_k(t)$, we calculate the first few of them explicitly:
\begin{eqnarray}
a_1(t) &=& a \,e^{-\nu t}, \label{a1t} \\
a_2(t) &=& -i a^2 \, \nu^{-1} \left[-3e^{-2\nu t} +3e^{-4\nu t} \right], \label{a2t} \\
a_3(t) &=& - a^3 \,\nu^{-2} \left[9e^{-3\nu t} -\frac{27}{2} e^{-5\nu t} + \frac{9}{2}
 e^{-9\nu t}\right],\label{a3t}\\
a_4(t) &=& i a^4 \, \nu^{-3} \left[-27e^{-4 \nu t} + 54e^{-6 \nu t} -\frac{27}{2}e^{-8\nu t} -18 e^{-10 \nu t}+
\frac{9}{2}e^{-16 \nu t} \right], \label{a4t}\\
a_5(t) &=& a^5 \, \nu^{-4}\left[81 e^{-5 \nu t}-\frac{405}{2} e^{-7 \nu t}+\frac{405}{4}
e^{-9 \nu t}+\frac{135}{2} e^{-11 \nu t}
-\frac{135}{4} e^{-13 \nu t} \right.
\nonumber\\
 && \, \left.-\frac{135}{8} e^{-17 \nu t}+\frac{27}{8} e^{-25 \nu t}\right ], \nonumber\\
a_6(t) &=& -ia^6 \, \nu^{-5}\left[-243 e^{-6 \nu t}+729 e^{-8 \nu t}-\frac{2187}{4}
 e^{-10 \nu t}-\frac{729}{4} e^{-12 \nu t} \right.
\nonumber\\
 &&
 \, \left. +243 e^{-14 \nu t}+\frac{81}{2} e^{-18 \nu t}
 -\frac{243}{8} e^{-20 \nu t}-\frac{243}{20} e^{-26 \nu t}+\frac{81}{40} e^{-36 \nu t}
\right].\nonumber
\end{eqnarray}

\vspace{.1in}
The following lemma summarizes the pattern exhibited by  $a_k(t)$'s.
\begin{lemma}\label{pattern}
For any $t>0$,
\begin{equation}\label{pai}
a_1(t) =a\,b_1(t), \quad a_2(t) = ia^2 \,b_2(t), \quad a_3(t) = -a^3 \,b_3(t),\quad a_4(t) = -i a^4 \,b_4(t)
\end{equation}
and more generally, for $k=4n+j$ with $n=0,1,2,\cdots$ and $j=1,2,3,4$,
\begin{equation}\label{pa}
a_{k}(t)=a_{4n+j}(t) = i^{j-1}\, a^{4n+j}\, b_{4n+j}(t),
\end{equation}
where $b_{4n+j}(t)>0$ for any $t>0$.
\end{lemma}

\noindent{\bf Remark}. A special consequence of this lemma is that all terms in the summation in (\ref{iter})
have the same sign and thus
\begin{equation} \label{abso}
|a_k(t)| = 3k e^{-\nu k^2 t} \, \int_0^t e^{\nu k^2 \tau} \sum_{k_1+k_2=k} |a_{k_1}(\tau)|\,|a_{k_2}(\tau)|\,d\tau.
\end{equation}
{\it Proof of Lemma \ref{pattern}}. \,\,(\ref{pa}) can be shown through induction. For $n=0$, (\ref{pa})
is just (\ref{pai}). By (\ref{iter}), $a_1(t) = a \,e^{-\nu t}$ and
$$
a_2(t) = 6i\,a^2 e^{-4\nu t} \int_0^t e^{4\nu \tau} \, b_1^2(\tau) \, d\tau =ia^2 b_2(t),
$$
where $b_2(t) = 6 e^{-4\nu t} \int_0^t e^{4\nu \tau} \, b_1^2(\tau) \, d\tau >0$. Similarly, $a_3(t) = -a^3 \,b_3(t)$
and $a_4(t) = -i a^4 \,b_4(t)$ for some $b_3(t)>0$ and $b_4(t)>0$.

\vspace{.1in}
We now consider the general case.  Without loss of generality, we prove (\ref{pa}) with $k=4n+1$. Assume (\ref{pa})
is true for all $k < 4n+1$. By (\ref{iter}),
$$
a_k(t) = 3ik\, e^{-\nu k^2\,t} \int_0^t e^{\nu k^2 \,\tau} \sum_{k_1+k_2=k} a_{k_1}(\tau)\, a_{k_2}(\tau) \, d\tau.
$$
Noticing that $a_{k_1}(\tau)\, a_{k_2}(\tau)$ with $k_1+k_2=4n+1$ assumes two forms
$$
a_{4n_1}(\tau)\, a_{4n_2+1}(\tau) \quad \mbox{and}\quad a_{4n_1+2}(\tau)\, a_{4n_2-1}(\tau)
$$
where $n_1\ge 0$, $n_2\ge 0$ and $n_1+n_2=n$, we conclude by the inductive assumptions
 that $a_{k_1}(\tau)\, a_{k_2}(\tau)$ must be of the form $-i\,a^k b_{k1,k_2}(\tau)$
for some positive function $b_{k1,k_2}(\tau)>0$. Therefore,
$$
a_k(t) =a_{4n+1}(t) = a^k \, b_{k}(t)
$$
with
$$
b_{k} (t) =3k\, e^{-\nu k^2\,t} \int_0^t e^{\nu k^2 \,\tau} \sum_{k_1+k_2=k} b_{k1,k_2}(\tau)\,d\tau >0
\quad\mbox{for any $t>0$}.
$$
This completes the proof of Lemma \ref{pattern}.

\vspace{.1in}
\noindent {\it Proof of Theorem \ref{bup}}. \,\, Without loss of generality, we set $\nu=1$. Assume
\begin{equation}\label{t0}
T\ge T_0 \equiv \sum_{k=2}^\infty \frac{1}{k^2} \ln\frac{3k-3}{2k-3}
\end{equation}
and choose $a$ such that
$$
A\equiv a\,e^{- T} \ge 1
$$
We prove by induction that
\begin{equation}\label{ub}
|a_k(T)| \ge A^k \quad\mbox{for}\quad k=1,2,3,\cdots
\end{equation}
which, in particular, yields (\ref{l2b}). Obviously, for any $0\leq
t\leq T$,
$$
|a_1 (t)|\geq |a_1(T)| = a \, e^{- T} = A \ge 1.
$$
To prove (\ref{ub}) for $k\ge 2$, we recall (\ref{abso}), namely
$$
|a_k(t)| = 3k e^{- k^2 t} \int_0^t e^{k^2\tau} \, \sum_{k_1+k_2=k} |a_{k_1}(\tau)|\,|a_{k_2}(\tau)|\, d\tau.
$$
Therefore, for $T\geq t\ge t_2 \equiv \frac14\,\ln3$,
$$
|a_2(t)| = 6 e^{-4 t} \int_0^t e^{4\tau} a_1^2(\tau) \,d\tau = \frac{3}{2} A^2 (1-e^{-4 t}) \ge A^2.
$$
For $k=3$, if $T\geq t\ge t_3 \equiv t_2 + \frac19 \ln\,2$,
\begin{eqnarray}
|a_3(t)| &=& 9 e^{-9 t} \int_0^t e^{9\tau} 2|a_1(\tau)|\,|a_2(\tau)| \,d\tau
\nonumber \\
&\ge& 9 e^{-9 t} \int_{t_2}^t e^{9\tau} 2|a_1(\tau)\,a_2(\tau)| \,d\tau \nonumber \\
&\ge& 2 A^3 \,(1-e^{-9(t-t_2)}) \ge A^3. \nonumber
\end{eqnarray}
More generally, for any $\geq t\ge t_k=t_{k-1} + \frac{1}{k^2} \,\ln
\frac{3k-3}{2k-3}$,
\begin{eqnarray}
|a_k(t)| &=& 3 k e^{-k^2 t} \,\int_0^t e^{k^2\tau}\, (|a_1(\tau)|\,|a_{k-1}(\tau)| + |a_2(\tau)|\,|a_{k-2}(\tau)|
\nonumber \\
&& + \cdots + |a_{k-2}(\tau)||a_2(\tau)| + |a_{k-1}(\tau)|\,|a_1(\tau)|) \,d\tau
\nonumber \\
& \ge & 3 k e^{-k^2 t} \,\int_{t_{k-1}}^t e^{k^2\tau}\, (|a_1(\tau)|\,|a_{k-1}(\tau)| + |a_2(\tau)|\,|a_{k-2}(\tau)|
\nonumber \\
&& + \cdots + |a_{k-2}(\tau)||a_2(\tau)| + |a_{k-1}(\tau)|\,|a_1(\tau)|) \,d\tau
\nonumber \\
&\ge& \frac{3k(k-1)}{k^2} (1-e^{-\nu\,k^2 (t-t_{k-1})}) A^k  \ge A^k. \nonumber
\end{eqnarray}
If $T\ge T_0$ as defined in (\ref{t0}), then $t_k < T$ for any integer $k\ge 1$ and thus
$$
|a_k(T)| \ge A^k.
$$
This completes the proof of Theorem \ref{bup}.

\vspace{.1in}
We state and prove a few specific properties for $a_k(t)$.

\begin{prop}\label{p} Assume $u_0$ is given by (\ref{ini}).
For each $k\geq 1$, $a_k(t)$ is of the form
\begin{equation} \label{formula}
a_k(t)=\sum_{m=k}^{k^2}\alpha_{k,m}e^{-m\nu t}
\end{equation}
where the complex-valued coefficients $\alpha_{k,m}$ satisfy
\begin{eqnarray}
&& \sum_{m=k}^{k^2} \alpha_{k,m} =0 \quad\mbox{for $k\ge 2$},\label{ak2}\\
&& \alpha_{k,m} = \frac{3ik}{k^2-m} \sum_{k_1+k_2=k}\sum_{m_1+m_2=m} \alpha_{k_1,m_1}\alpha_{k_2,m_2}
 \quad\mbox{for $k\le m <k^2$}. \label{coe}
\end{eqnarray}
The indices $k_1$, $k_2$, $m_1$ and $m_2$ in the summation above obey
$$
1\le k_1\le k-1,\,\, 1\le k_2\le k-1,\,\, k_1\le m_1\le k_1^2 \quad \mbox{and}\quad k_2 \le m_2 \le k_2^2.
$$
\end{prop}
{\it Proof}. (\ref{ak2}) is a consequence of the fact that $a_k(0)=0$ for $k\ge 2$. (\ref{formula})
follows from a simple induction. Obviously, $a_1(t)=a\,e^{-\nu t}$.
Fix $k$ and assume (\ref{formula}) is valid for all integers up to $k$. Then, for
$k_1\geq 1, k_2\geq 1, k_1+k_2=k+1, k_1\le m_1\leq k_1^2$ and $k_2\le m_2\leq k_2^2$,
\begin{eqnarray}
& a_{k+1}(t) &= 3i(k+1)\,\sum_{k_1+k_2=k+1}\sum_{m_1,m_2}\alpha_{k_1,m_1}\alpha_{k_2,m_2}
 e^{-\nu (k+1)^2\,t} \int_0^t e^{\nu ((k+1)^2-(m_1+m_2))\,\tau} \, d\tau \nonumber\\
 & &=\sum_{k_1+k_2=k+1} \sum_{m_1,m_2}\frac{3i(k+1)\,\alpha_{k_1,m_1}\alpha_{k_2,m_2}}{\nu ((k+1)^2-(m_1+m_2))}
\left(e^{-\nu(m_1+m_2)t}-e^{-\nu (k+1)^2 t}\right).\nonumber
\end{eqnarray}
Since $m_1+m_2\leq k_1^2+k_2^2\leq (k_1+k_2)^2=(k+1)^2$, this proves (\ref{formula}) with (\ref{coe}).

\vspace{.1in}
\begin{prop}\label{prop2}
Assume that $u_0$ is given by (\ref{ini}).
\begin{enumerate}
\item[1)] Let $k\ge 1$ be an integer.  Then
\begin{equation}\label{kk}
\alpha_{k,k}=\left(\frac{3i}{\nu}\right)^{k-1}a^k  \quad \mbox{and} \quad \alpha_{k,k+2}=-\frac{k}{2}\,\alpha_{k,k};
\end{equation}
\item[2)] Let $k\ge 1$ be an integer.  Then, for $n=1,3,5,\cdots,$
$$
\alpha_{k,k+n} =0;
$$
\item[3)] Let $k\ge 1$ be an integer and let $k^2>m>U(k)\equiv k^2-2k+2$. Then
\begin{equation}\label{km0}
\alpha_{k,m} =0.
\end{equation}
\end{enumerate}
\end{prop}
{\it Proof}. Letting $m_1=k_1$ and $m_2=k_2$ in (\ref{coe}), we find
\begin{eqnarray}
\alpha_{k,k}=\sum_{k_1+k_2=k}\alpha_{k_1,k_1}\alpha_{k_2,k_2} \frac{3ik}{\nu (k^2-k)}
            =\frac{3i}{\nu (k-1)}\sum_{k_1+k_2=k}\alpha_{k_1,k_1}\alpha_{k-k_1,k-k_1}. \nonumber
\end{eqnarray}
A simple induction allows us to obtain the expression for $\alpha_{k,k}$.
To show $\alpha_{k,k+2}=-\frac{k}{2}\,\alpha_{k,k}$, we set $m=k+2$ in (\ref{coe}) to get
\begin{eqnarray}
\alpha_{k,k+2} &=& \frac{3ik}{\nu(k^2-k-2)}(\alpha_{1,1}\, \alpha_{k-1,k+1} + \alpha_{2,2}\,
\alpha_{k-2,k} + \alpha_{2,4}\,\alpha_{k-2,k-2} \nonumber\\
&& +\cdots\,+ \alpha_{k-2,k-2}\,\alpha_{2,4} + \alpha_{k-2,k}\, \alpha_{2,2} + \alpha_{1,1}\,
\alpha_{k-1,k+1}). \label{2k2}
\end{eqnarray}
Inserting the inductive assumptions such as
$$
\alpha_{k-1,k+1}=-\frac{k-1}{2}\, \alpha_{k-1,k-1},\quad \alpha_{k-2,k}=-\frac{k-2}{2} \alpha_{k-2,k-2},
 \quad \alpha_{2,4}=-\alpha_{2,2}
$$
in (\ref{2k2}), we obtain
\begin{eqnarray}
\alpha_{k,k+2} &=& \frac{3ik}{\nu(k^2-k-2)}\left[-\frac{k}{2} \sum_{k_1=1}^{k-1} \alpha_{k_1,k_1}\,
\alpha_{k-k_1,k-k_1} + \alpha_{1,1}\,\alpha_{k-1,k-1}\right] \nonumber\\
&=& -\frac{k}{2} \frac{k^2-k}{k^2-k-2} \frac{3ik}{\nu(k^2-k)} \sum_{k_1=1}^{k-1} \alpha_{k_1,k_1}\,
\alpha_{k-k_1,k-k_1}
\nonumber\\
&& \qquad + \frac{3ik}{\nu(k^2-k-2)} \alpha_{1,1}\,\alpha_{k-1,k-1} \nonumber\\
&=& -\frac{k}{2} \frac{k^2-k}{k^2-k-2} \alpha_{k,k} -\frac{k}{2} \frac{-2}{k^2-k-2} \alpha_{k,k}
 =-\frac{k}{2}\,\alpha_{k,k}.\nonumber
\end{eqnarray}
To show $\alpha_{k,k+1}=0$, we set $m=k+1$ to obtain
$$ \alpha_{k,k+1}=\frac{3ik}{\nu(k^2-(k+1))} \,\left(\alpha_{1,1}\,\alpha_{k-1,k} + \alpha_{2,2}
\,\alpha_{k-2,k-1} +\cdots + \alpha_{k-1,k}\,\alpha_{1,1} \right),
$$
which can be seen to be zero after inserting the inductive assumptions.

\vspace{.1in}
To prove (\ref{km0}), it suffices to notice in (\ref{coe}) that the second summation is over $m_1+m_2=m$ with $k_1\le m_1\le k_1^2$ and $k_2\le m_2\le k_2^2$. Thus, $m=m_1+m_2 \le k_1^2 + k_2^2 =(k_1+k_2)^2 -2k_1k_2 \le k^2-2(k-1)$ and $\alpha_{k,m}$ with $U(k)<m<k^2$ is equal to zero. This completes the proof of Proposition \ref{prop2}.

\vspace{.15in}
\noindent{\bf \ref{sec:blow}.2\,\, Local well-posedness}

\vspace{.1in}
This subsection establishes the following two major results.

\begin{thm}\label{local1} Consider (\ref{kdvb}) with $\gamma>\frac12$.
Let $s>\frac12$. Assume $u_0\in H^s({\mathbf T})$ has the form
\begin{equation}\label{inis}
u_0(x) =\sum_{k=1}^\infty \,a_{0k}\, e^{ikx}.
\end{equation}
Then there exists $T=T(\|u_0\|_{H^s})$ such that (\ref{kdvb}) with the initial data $u_0$ has
a unique solution $u\in C([0,T); H^s)\cap L^2([0,T);\mathring{H}^{s+\gamma})$ that assumes the form
$$
u(x,t) =\sum_{k=1}^\infty \, a_k(t)\, e^{ikx}.
$$
\end{thm}

\vspace{.1in}
In the case when $\gamma\ge 1$, we can actually show that no finite-time singularity is possible if we know that the $L^2$-norm is bounded {\it a priori}. In fact, the following theorem states that the $L^2$-norm controls all higher-order derivatives.

\begin{thm} \label{nu}
Let $T>0$ and let $u$ be a weak solution of (\ref{kdvb}) with $\gamma\ge 1$ on the time interval $[0,T]$. If we know {\it a priori} that $u\in L^\infty([0,T]; L^2) \cap L^2([0,T]; \mathring{H}^\gamma)$, namely
\begin{equation}\label{M0}
M_0\equiv \sup_{t\in[0,T]} \|u(\cdot,t)\|^2_{L^2} +\nu \int_0^T \|\Lambda^\gamma u(\cdot,t)\|^2_{L^2} \,dt <\infty,
\end{equation}
then, for any integer $k>0$,
$$
M_k \equiv \sup_{t\in[0,T]} \|u^{(k)}(\cdot,t)\|^2_{L^2} + \nu \int_0^T \|\Lambda^{k+\gamma} u(\cdot,t)\|^2_{L^2} \,dt <\infty.
$$
where $\Lambda =(-\Delta)^\frac12$ and $u^{(k)}$ denotes any partial derivative of order $k$.
\end{thm}

\vspace{.1in}
We first prove Theorem \ref{local1}.

\noindent {\it Proof of Theorem \ref{local1}}.\quad The existence of such a solution follows from the
Galerkin approximation. Let $N\ge 1$ and denote by $P_N$ the projection on the subspace $\{e^{ix},
e^{2ix},\cdots, e^{iNx}\}$. Let
$$
u^N(x,t) = \sum_{k=1}^N a^N_k(t) \,e^{ikx}
$$
where $a_k(t)$ satisfies
\begin{eqnarray}
&& \frac{d}{dt} \,a^N_k(t) = 3ik \, \sum_{k_1+k_2=k} a^N_{k_1}(t) \,a^N_{k_2}(t) + i \alpha\,k^3\,a^N_k(t)
-\nu k^{2 \gamma} \,a^N_k(t), \nonumber \\
&& a^N_k(0) = a^N_{0k}\equiv a_{0k}. \label{ode}
\end{eqnarray}
Here $1\le k_1 \le N$ and $1\le k_2\le N$. From the theory of ordinary differential equations, we know that
 (\ref{ode}) has a unique local solution $a^N_k(t)$ on $[0,T]$. We derive some {\it a priori} bounds for $u^N(x,t)$.
 Clearly, $u^N(x,t)$ solves
$$
\partial_t u^N =6 P_N(u^N u^N_x) +\alpha u^N_{xxx}- \nu (-\Delta)^\gamma u^N, \quad u^N(x,0) = P_N \,u_0.
$$
We now show that
\begin{equation}\label{gg}
\frac{d}{dt} \|u^N\|_{H^s}^2 + \nu \|u^N\|^{2}_{H^{s+\gamma}} \le C(\nu,s) \|u^N\|_{H^s}^\frac{6\gamma-2}{2\gamma-1}.
\end{equation}
It follows from the equation
$$
\frac{d}{dt} a^N_k(t) + \nu k^{2\gamma} a^N_k(t) - i\alpha k^3\,a^N_k(t) = 3 ik \sum_{k_1+k_2=k}
 a^N_{k_1}(t)\,a^N_{k_2}(t)
$$
that, after omitting the upper index $N$ for notational convenience,
$$
\frac{d}{dt} \sum_{k=1}^N k^{2s} \, |a_k(t)|^2 = -2 \nu \sum_{k=1}^N k^{2(s+\gamma)} \,
 |a_k(t)|^2 -6 \sum_{k=1}^N k^{2s+1} \, {\mathcal I}\left(\bar{a}_k\,\sum_{k_1+k_2=k} a_{k_1}\,a_{k_2}\right),
$$
where ${\mathcal I}$ denotes the imaginary part. To bound the nonlinear term on the right
(denoted by $J$), we first notice that the summation over $k_1+k_2=k$ is less than twice
 the summation over $k_1+k_2=k$ with  $k_1\le k_2$ and $2k_2\ge k$. Thus,
\begin{eqnarray}
J &\le& 6 \sum_{k=1}^N k^{2s+1} \, |a_k|\,\sum_{k_1+k_2=k} |a_{k_1}|\,|a_{k_2}| \nonumber \\
&\le & 12 \sum_{k=1}^N k^{s+\frac12} |a_k|\,\sum_{k/2\le k_2 \le k} (2k_2)^{s+\frac12} |a_{k_1}|\,|a_{k_2}|. \nonumber
\end{eqnarray}
Applying H\"{o}lder's inequality and Young's inequality for series, we have
\begin{eqnarray}
J &\le& 12\,\left[\sum_{k=1}^N k^{2s+1} |a_k|^2\right]^\frac12\, \left[\sum_{k=1}^N
 \left(\sum_{k/2\le k_2 \le k} (2k_2)^{s+\frac12} |a_{k_1}|\,|a_{k_2}|\right)^2 \right]^\frac12
  \nonumber\\
&\le& 12 \left[\sum_{k=1}^N k^{2s+1} |a_k|^2\right]^\frac12\, \left[\sum_{k_2=1}^N k_2^{2s+1}
 |a_{k_2}|^2\right]^\frac12\, \sum_{k_1=1}^N |a_{k_1}|.  \nonumber \\
&\le& 12 \sum_{k=1}^N k^{2s+1} |a_k|^2 \, \left[\sum_{k_1=1}^N |k_1|^{2s} |a_{k_1}| \right]^{\frac12}\,
 \left[\sum_{k_1=1}^N k_1^{-2s}\right]^{\frac12} \nonumber \\
&\le& C(s) \|u^N\|^2_{\mathring{H}^{s+\frac12}}\, \|u^N\|_{H^s}.
\end{eqnarray}
Thus, we get
\begin{equation}\label{ddt}
\frac{d}{dt} \|u^N\|_{H^s}^2 + 2 \nu \|u^N\|^{2}_{\mathring{H}^{s+\gamma}} \le C(s)
 \|u^N\|^2_{\mathring{H}^{s+\frac12}}\, \|u^N\|_{H^s}
\end{equation}
By H\"{o}lder's inequality
$$
\|u^N\|_{\mathring{H}^{s+\frac12}} \le \|u^N\|^{\frac1{2\gamma}}_{\mathring{H}^{s+\gamma}}\,
\|u^N\|_{H^s}^{1-\frac1{2\gamma}},
$$
we have
\begin{equation}\label{ddt1}
J\le C(s) \|u^N\|^{\frac1{\gamma}}_{\mathring{H}^{s+\gamma}}\,\|u^N\|_{H^s}^{3-\frac1{\gamma}}
\le \nu \|u^N\|^{2}_{\mathring{H}^{s+\gamma}}\, + C(\nu,s)\, \|u^N\|_{H^s}^\frac{6\gamma-2}{2\gamma-1}.
\end{equation}
(\ref{ddt}) and (\ref{ddt1}) yield (\ref{gg}). With these bounds at our disposal, the existence of
a solution $u$ of the form (\ref{form}) is then obtained as a limit of $u^N$ as $N\to \infty$.

\vspace{.1in}
We now turn to the uniqueness. Assume (\ref{kdvb}) has two solutions $u_1$ and $u_2$ satisfying
$$
u_1,\,\, u_2 \in C([0,T); H^s)\cap L^2([0,T);\mathring{H}^{s+\gamma}).
$$
Then their difference $w=u_1-u_2$ satisfies
$$
w_t + \nu (-\Delta)^\gamma w + \alpha w_{xxx} = 6 w u_{1x} + 6 u_2 w_x.
$$
Applying the same procedure as in the derivation of (\ref{ddt}), we find that, for $s>\frac12$,
$$
\frac{d}{dt} \|w\|_{H^s}^2 + 2 \nu \|w\|^{2}_{\mathring{H}^{s+\gamma}} \le C(s)
\,\|w\|^2_{H^s}(\|u_1\|_{\mathring{H}^1} + \|u_2\|_{\mathring{H}^1}).
$$
The fact that $u_1, u_2\in L^2([0,T);\mathring{H}^{s+\gamma})$ with $s+\gamma>1$ and an application of
Gronwall's inequality yields the uniqueness.
This completes the proof of Theorem \ref{local1}.

\vspace{.1in}
\noindent{\it Proof of Theorem \ref{nu}}. \,\, We start with the case $k=1$.  It is easy to verify that
\begin{equation}\label{base}
\frac{d}{dt} \|u_x(\cdot,t)\|_{L^2}^2 + 2 \kappa \|\Lambda^\gamma u_x \|_{L^2}^2
= I_1 + I_2,
\end{equation}
where
\begin{eqnarray}
I_1 &=& 2 \int |u_x|^2 {\mathcal R}(u_x)\,dx, \nonumber\\
I_2 &=& 2 \int {\mathcal R}(u\, \overline{u}_x\, u_{xx})\, dx. \nonumber
\end{eqnarray}
Here ${\mathcal R}$ denotes the real part. By the Gagliardo-Nirenberg type equalities,
\begin{eqnarray}
|I_1| &\le & 2 \|u_x\|_{L^2}^2 \|u_x\|_{L^\infty} \nonumber\\
& \le & C \|u\|_{L^2}^{\gamma_1} \, \|u_x\|_{L^2}^2\,  \|\Lambda^{1+\gamma}u\|_{L^2}^{1-\gamma_1}, \nonumber
\end{eqnarray}
\begin{eqnarray}
|I_2| &\le & C \|u\|_{L^\infty} \,\|u_x\|_{L^2} \, \|u_{xx}\|_{L^2} \nonumber\\
& \le & C \|u\|_{L^2}^{\frac12} \, \|u_x\|_{L^2}^{\frac32}\, \|u_{xx}\|_{L^2} \nonumber\\
& \le & C \|u\|_{L^2}^{\frac12 + \gamma_2} \, \|u_x\|_{L^2}^{\frac32}\, \|\Lambda^{1+\gamma} u\|_{L^2}^{1-\gamma_2}
\nonumber
\end{eqnarray}
where
$$
\gamma_1= \frac{2\gamma-1}{2\gamma+2} \quad \mbox{and} \quad \gamma_2= \frac{\gamma-1}{\gamma+1}.
$$
By Young's inequality,
\begin{eqnarray}
|I_1| &\le & \frac{\nu}{2} \|\Lambda^{1+\gamma}u\|_{L^2}^2 + C \nu^{-\frac{1-\gamma_1}{1+\gamma_1}}\,
\|u\|_{L^2}^{\frac{2\gamma_1}{1+\gamma_1}}\, \|u_x\|_{L^2}^{\frac{4}{1+\gamma_1}}, \nonumber \\
|I_2| &\le & \frac{\nu}{2} \|\Lambda^{1+\gamma}u\|_{L^2}^2 + C \nu^{-\frac{1-\gamma_2}{1+\gamma_2}}\,
 \|u\|_{L^2}^{\frac{1+2\gamma_2}{1+\gamma_2}}\, \|u_x\|_{L^2}^{\frac{3}{1+\gamma_2}}.\nonumber
\end{eqnarray}
Inserting these inequalities in (\ref{base}) and integrating with respect to $t$ yields
\begin{eqnarray}
&& \sup_{t\in[0,T]} \|u_x(\cdot,t)\|^2_{L^2}  + \nu \int_0^T \|\Lambda^{1+\gamma} u\|^2_{L^2} \,dt \nonumber \\
&& \qquad \qquad \qquad
\le C(\nu) M_0^{\frac{\gamma_1}{1+\gamma_1}}\, \int_0^T \|u_x\|_{L^2}^{\frac{4}{1+\gamma_1}}\,dt
+ C(\nu)\, M_0^{\frac{1+2\gamma_2}{2+2\gamma_2}}\, \int_0^T \|u_x\|_{L^2}^{\frac{3}{1+\gamma_2}}\,dt, \nonumber
\end{eqnarray}
where $M_0$ is specified in (\ref{M0}). By (\ref{M0}) and the Gagliardo-Nirenberg type inequality
$$
\|u_x\|_{L^2} \le C \|u\|_{L^2}^{1-\frac1\gamma}\, \|\Lambda^\gamma u\|_{L^2}^{\frac1\gamma},
$$
we have
$$
\int_0^T \|u_x\|^{2\gamma}_{L^2}\,dt \le C M_0^\gamma.
$$
Therefore,
\begin{eqnarray}
&& \sup_{t\in[0,T]} \|u_x(\cdot,t)\|^2_{L^2}  + \nu \int_0^T
\|\Lambda^{1+\gamma} u\|^2_{L^2} \,dt \nonumber \\
&& \qquad \qquad \qquad
\le C(\nu) M_0^{\frac{4\gamma^2+3\gamma-1}{4\gamma+1}}\,\sup_{t\in[0,T]}
\|u_x(\cdot,t)\|_{L^2}^{\frac{-8\gamma^2+6\gamma+8}{4\gamma+1}} \nonumber \\
&& \quad \qquad \qquad \qquad
+\, C(\nu) M_0^{\frac{4\gamma^2+3\gamma-1}{4\gamma}}\,\sup_{t\in[0,T]}
\|u_x(\cdot,t)\|_{L^2}^{\frac{-4\gamma^2+3\gamma+3}{2\gamma}}. \label{good1}
\end{eqnarray}
When $\gamma > \frac34$, $4\gamma^2+\gamma-3 >0$ and consequently
$$
\frac{-8\gamma^2+6\gamma+8}{4\gamma+1} <2 \quad \mbox{and}\quad \frac{-4\gamma^2+3\gamma+3}{2\gamma}<2.
$$
(\ref{good1}) then implies that
$$
\sup_{t\in[0,T]} \|u_x(\cdot,t)\|^2_{L^2}  + \nu \int_0^T \|\Lambda^{1+\gamma} u\|^2_{L^2} \,dt \le M_1,
$$
where $M_1$ is a constant depending on $\gamma$, $\nu$ and $M_0$ alone. $L^2$-bounds
 for higher-order derivatives can be obtained through iteration. This completes the proof of Theorem \ref{nu}.

\vspace{.2in}
\section{Global solutions of the complex KdV-Burgers equation}
\setcounter{equation}{0}
\label{series}

We consider the initial-value problem for the complex KdV-Burgers equation
\begin{equation} \label{kdvp}
\left\{
\begin{array}{l}
\displaystyle u_t -6 uu_x + \alpha u_{xxx} -\nu u_{xx} =0,\quad x\in {\mathbf T}, \,t>0, \\
\displaystyle u(x,0) =u_0(x), \quad x\in {\mathbf T}
\end{array}
\right.
\end{equation}
and study the global regularity of its solutions of the form
\begin{equation}\label{ss}
u(x,t) = \sum_{k=1}^\infty a_k(t) e^{ikx}.
\end{equation}
Here $\alpha\ge 0$ and $\nu\ge 0$ and (\ref{kdvp}) includes the complex Burgers and complex
 KdV equations as special cases. Two major results are established. Theorem \ref{series2}
 presents a general conditional global regularity result and  Theorem \ref{sing} asserts
 the global regularity of (\ref{ss}) for a special case.

\vspace{.1in}
Assume the initial data $u_0$ is of the form
\begin{equation}\label{ind}
u_0(x) = \sum_{k=1}^\infty \, a_{0k}\,e^{ikx}
\end{equation}
and is in $H^s$ with $s>\frac12$. According to Theorem \ref{local1}, (\ref{kdvp})
 has a unique local solution $u\in C([0,T); H^s)$ of the form (\ref{ss}) for some $T>0$.
To study the global regularity of (\ref{ss}), we explore the structure of $a_k(t)$
and obtain the following two propositions.
\begin{prop} \label{for}
If (\ref{ss}) solves (\ref{kdvp}), then
$a_k(t)$ can be written as
\begin{equation} \label{akform}
a_k(t) = \sum_{ k\le h \le k^2,\, k\le l \le k^3 }
a_{k,\,h,\,l} \, e^{-(\nu h - \alpha i l) t}
\end{equation}
where $a_{k,\,h,\,l}$ consists of a finite number of terms of the form
\begin{equation} \label{ak1form}
C(\alpha,\nu,k,h,l, j_1,\cdots, j_k)\, a_{01}^{j_1} \,a_{02}^{j_2} \, \cdots \, a_{0k}^{j_k}
\end{equation}
with $j_1$, $j_2$ ,$\cdots$, $j_k$ being nonnegative integers and satisfying
\begin{equation} \label{jcon}
 j_1 + 2 j_2 + \cdots + kj_k =k.
\end{equation}
\end{prop}

\begin{prop} \label{for2}
Let $k\ge 1$ be an integer. Let $U(k)=k^2-2k+2$ and $V(k)=k^3-3k^2+3k$.  The
coefficients\, $a_{k,h,l}$ in (\ref{akform}) have the following properties
\begin{enumerate}
\item[(1)] For $k\le h<k^2$ and $k\le l < k^3$,
\begin{equation} \label{ty1}
a_{k,h,l} = \frac{3ik}{\nu(k^2-h)-i\alpha(k^3-l)} \,\sum_{k_1+k_2=k}\sum_{h_1+h_2=h}
 \sum_{1_1+1_2=l} \alpha_{k_1,h_1,l_1}\alpha_{k_2,h_2,l_2}
\end{equation}
\item[(2)] For $h=k^2$ and $l=k^3$,
\begin{equation} \label{ty2}
a_{k,k^2,k^3} = a_k(0)- \sum_{k\le h<k^2}\sum_{k\le l<k^3}  a_{k,h,l}
\end{equation}
\item[(3)] For $U(k)<h<k^2$ or $V(k)<l<k^3$,
\begin{equation} \label{ty3}
a_{k,h,l} =0.
\end{equation}
\end{enumerate}
\end{prop}

\noindent {\it Proof of Proposition \ref{for}}.\,\,If (\ref{ss}) solves (\ref{kdvp}) ,
then $a_k(t)$ solves the
ordinary differential equation
$$
\frac{d}{dt} a_k(t) + (\nu k^2 -\alpha i k^3) a_k(t) -3ik \sum_{k_1+k_2=k}
a_{k_1}(t) \,a_{k_2}(t) =0.
$$
The equivalent integral form is given by
\begin{equation} \label{mit}
a_k(t) = e^{-(\nu k^2 -\alpha i k^3)t} \left[a_{0k} +3 ik\int_0^t e^{(\nu k^2 -\alpha i k^3)\tau}
\sum_{k_1+k_2=k} a_{k_1}(\tau) \,a_{k_2}(\tau)\,d\tau \right].
\end{equation}
It is easy to show through an inductive process that $a_k$ is of the form (\ref{akform}).
 In addition, for $k\le h<k^2$ and $k\le l<k^3$, the term  in (\ref{ak1form}) with fixed $j_1$,
 $j_2$, $\cdots$, $j_k$ satisfying
$$
j_1+2j_2 +\cdots + kj_k=k
$$
can be expressed as
\begin{eqnarray}
&& C(\alpha,\nu,k,h,l, j_1,\cdots, j_k)\, a_{01}^{j_1} \,a_{02}^{j_2} \, \cdots \, a_{0k}^{j_k}\nonumber\\
&& \qquad
=\frac{3ik}{\nu(k^2-h)-i\alpha(k^3-l)}\sum_{m_1+n_1=j_1}\cdots\sum_{m_k+n_k=j_k}
C(\alpha,\nu,k_1,h_1,l_1, m_1,\cdots, m_{k_1})\nonumber\\
&& \qquad\quad \times \, C(\alpha,\nu,k_2,h_2,l_2, n_1,\cdots, n_{k_2})
\,a_{01}^{m_1+n_1} \,a_{02}^{m_2+n_2} \, \cdots \, a_{0k}^{m_k+n_k} \label{khl}
\end{eqnarray}
where the indices satisfy
\begin{eqnarray}
&& 1\le k_1\le k, \quad  1\le k_2\le k,\quad k_1+k_2=k, \nonumber\\
&& k_1\le h_1\le k_1^2, \quad k_2\le h_2 \le k_2^2, \quad h_1+h_2=h, \nonumber\\
&& k_1\le l_1 \le k_1^3,\quad k_2\le l_2 \le k_2^3, \quad l_1+l_2=l, \nonumber\\
&& m_1+n_1=j_1,\quad m_2+n_2 =j_2, \quad \cdots, \quad m_k+n_k =j_k. \nonumber \\
&& (m_r=0 \quad\mbox{for $r>k_1$ and}\quad n_r=0 \quad\mbox{for $r>k_2$}) \nonumber\\
&& m_1+2m_2 +\cdots +k_1m_{k_1}=k_1,\quad n_1+2n_2 +\cdots +k_2n_{k_2}=k_2.\nonumber
\end{eqnarray}
When $h=k^2$ and $l=k^3$,
\begin{equation}\label{kk2}
C(\alpha,\nu,k,k^2,k^3, j_1,j_2,\cdots, j_k)=\left\{
\begin{array}{l}
1 \quad \mbox{for}\quad (j_1,j_2,\cdots,j_k)=(0,0,\cdots,1), \\
-C(\alpha,\nu,k,h,l, j_1,j_2,\cdots, j_k)\quad \mbox{otherwise}
\end{array}
\right.
\end{equation}
for some $h<k^2$ and $l<k^3$. To illustrate these formulas, we list $a_k$ for $k=1,2,3$,
\begin{eqnarray}
&& a_1(t) =a_{01}\, e^{-(\nu-i\alpha)t},
\nonumber \\
&& a_2(t) =\frac{6i}{2\nu-6\alpha i} a_{01}^2\, e^{-(2\nu-2\alpha i)t}+ \left[a_{02}
-\frac{6i}{2\nu-6\alpha i} a_{01}^2\right] e^{-(4\nu-8i\alpha)t}
, \nonumber \\
&& a_3(t) =\frac{108 a_{01}^3}{(2\nu-6\alpha i)(6\nu-24\alpha i)} e^{(-3\nu + 3\alpha i)t}\nonumber\\
&&\qquad\quad  + \left[\frac{18 ia_{01}a_{02}}{4\nu-18\alpha i}
- \frac{108 a_{01}^3}{(2\nu-6\alpha i)(4\nu-18\alpha i)}\right]\,e^{(-5\nu+9i\alpha)t}\nonumber\\
&& \qquad \quad + \left[a_{03}-\frac{18 ia_{01}a_{02}}{4\nu-18\alpha i}
+ \frac{108 a_{01}^3}{(2\nu-6\alpha i)(4\nu-18\alpha i)}
-\frac{108 a_{01}^3}{(2\nu-6\alpha i)(6\nu-24\alpha i)} \right]\nonumber\\
&& \qquad \qquad \quad\times \,e^{(-9\nu+27\alpha i)t}. \nonumber
\nonumber
\end{eqnarray}

\vspace{.1in}
\noindent{\it Proof of Proposition \ref{for2}}.\,\, (\ref{ty1}) follows from a simple induction. (\ref{ty2})
 is obtained by set $t=0$ in (\ref{akform}). To show (\ref{ty3}), we notice that the second summation in (\ref{ty1})
  is over $h_1+h_2=h$ with $k_1\le h_1\le k_1^2$ and $k_2\le h_2\le k_2^2$ while the third summation is
  over $l_1+l_2=l$ with  $k_1\le l_1\le k_1^3$ and $k_2\le l_2\le k_2^3$. Thus,
\begin{eqnarray}
&& h=h_1+h_2\le k_1^2 + k_2^2 =k^2-2k_1\,k_2\le k^2-2(k-1) =U(k), \nonumber\\
&& l=l_1+l_2 \le k_1^3 + k_2^3 =k^3 - 3k\,k_1\,k_2 \le k^3-3k(k-1) =V(k).\nonumber
\end{eqnarray}
That means, $a_{k,h,l}=0$ if $U(k)<h<k^2$ and $V(k)<l<k^3$.

\vspace{.1in}
\begin{thm} \label{series2}
Consider (\ref{kdvp}) with $\nu>0$. Assume  $u_0\in H^s({\mathbf T})$ with $s>\frac12$
can be represented in the form (\ref{ind}) with
\begin{equation} \label{l2c}
|a_{0k}| \le 1, \quad k=1,2,\cdots
\end{equation}
If we have the uniform bound
\begin{equation} \label{cc0}
|C(\alpha,\nu,k,h,l, j_1,\cdots, j_k)| \le C_0(\alpha,\nu)
\end{equation}
for all $k\ge 1$, $k\le h<k^2$, $k\le l<k^3$ and $(j_1,j_2,\cdots, j_k)$ satisfying (\ref{jcon}),
 then (\ref{kdvp}) has a unique global solution $u$ given by (\ref{ss}).
In addition, for any $s\ge 0$, there are $T_0>0$ and $\delta>0$ such that for any $t\ge T_0$,
\begin{equation}\label{hsb}
\|u(\cdot,t)\|_{H^s} < \frac{C(\alpha,\nu,s)}{1-e^{-\nu t}}\, e^{-\delta \nu k t}
\end{equation}
where  $C$ is a constant depending on $\alpha$, $\nu$ and $s$ only.
\end{thm}

We remark that the assumption in (\ref{cc0}) can be verified for the case when $a_{01}>0$ and $a_{02}=a_{03}=\cdots=0$.
We assume that $\nu$ and $\alpha$ satisfy $\nu^2+ 9\alpha^2\ge 36$ and show by induction that
$$
| C(\alpha,\nu,k,h,l, j_1,\cdots, j_k)| \le 1.
$$
Since $a_{02}=a_{03}=\cdots=0$, these coefficients are nonzero only if $j_1=k$ and $j_2=j_3=\cdots =j_k=0$.
For any $k\le h<k^2$ and $k\le l<k^3$, we have, according to (\ref{khl}),
\begin{eqnarray*}
&& | C(\alpha,\nu,k,h,l, j_1,\cdots, j_k)| \\
&& \qquad\qquad \le \left|\frac{3ik}{\nu(k^2-h)-i\alpha(k^3-l)}\right|\, \sum_{m_1+n_1=j_1}
|C(\alpha,\nu,k_1,h_1,l_1, m_1,\cdots, m_{k_1})|\\
&& \qquad \qquad \quad \times \, |C(\alpha,\nu,k_2,h_2,l_2, n_1,\cdots,  n_{k_2})|.
\end{eqnarray*}
For $j_1=k$, the number of terms in the summation $m_1+n_1=j_1$ is at most $k$. By the inductive assumption,
$$
| C(\alpha,\nu,k,h,l, j_1,\cdots, j_k)| \le \frac{3k^2}{\sqrt{\nu^2(k^2-h)^2 + \alpha^2(k^3-l)^2}}
$$
Applying (\ref{ty3}), $h \le U(k)\equiv k^2-2k+2$ and $l\le V(k)\equiv k^3-3k^2+3k$ and
 thus $| C(\alpha,\nu,k,h,l, j_1,\cdots, j_k)| \le 1$ by taking into account the assumption
 on $\nu$ and $\alpha$. When $h=k^2$ and $l=k^3$, the boundedness of the coefficient follows from (\ref{kk2}).

\vspace{.1in}
The proof of Theorem \ref{series2} involves a very classical problem in number theory, namely
the number of integer solutions $(j_1,j_2,\cdots,j_k)$ to the equation defined in (\ref{jcon})
for a given positive integer $k$. This problem is not as simple as it may look like.
An upper bound and an asymptotic approximation for the number of nonnegative solutions
are given by G.H. Hardy and S. Ramanujan \cite{HR}, as stated in the following lemma.

\begin{lemma}\label{hardy}
Let $k>0$ be an integer and let $N_k$ denote the number of nonnegative solutions
to the equation
$$
j_1 + 2 j_2 + \cdots + kj_k =k.
$$
Then, for some constant $C_1$,
$$
N_k  < \frac{C_1}{k} \, e^{2\sqrt{2\, k}}.
$$
In addition, $N_k$ has the following asymptotic behavior:
$$
N_k \quad \sim \quad \frac{1}{4\sqrt{3} k} e^{\pi \sqrt{\frac{2k}{3}}},
\quad \mbox{as $k\to \infty$}.
$$
\end{lemma}

\vspace{.1in}
\noindent{\it Proof of Theorem \ref{series2}}. Applying (\ref{l2c}) and (\ref{cc0}),
we obtain the following bound for $a_{k,\,h,\, l}$ in (\ref{akform})
$$
|a_{k,\,h,\, l}| \le C_0(\alpha,\nu)\, N_k
 \le \frac{C_2}{k} \,  e^{2\sqrt{2\, k}},
$$
where $C_2=C_0C_1$ and we have used Lemma \ref{hardy}. Therefore,
\begin{eqnarray}
|a_k(t)| &\le& \sum_{k\le h \le k^2} \sum_{k\le l \le k^3} |a_{k,h,l}|\,
e^{- \nu h t} \nonumber \\
&\le & C_2\,(k^2-1)\,  e^{2\sqrt{2}\,\sqrt{k}}\,
\frac{e^{-\nu k t}}{1-e^{-\nu t}}. \label{basis}
\end{eqnarray}
For any fixed $t>0$, we can choose $K=K(\nu)$ and $0<M=M(\nu)<1$ such that
$$
|a_k(t)| \le \frac{C_2}{1-e^{-\nu t}}\,M^k \quad\mbox{for $k\ge K$}.
$$
Therefore, $u$ represented by (\ref{ss}) converges for any $t>0$. In addition, $u(\cdot,t)\in H^s$
for any $s\ge 0$. To see the exponential decay of $\|u(\cdot,t)\|_{H^s}$ for large time, we
 choose $T_0=T_0(\nu,s)$ such that for any $t\ge T_0$ and $k\ge 1$
$$
(1+k^2)^s |a_k(t)|^2 \le C_2 \, M_1^k \, \frac{e^{-\delta\,\nu k t}}{1-e^{-\nu t}},
$$
where $M_1>0$ and $\delta>0$ are some constants. This bound then implies (\ref{hsb}).
This completes the proof of Theorem  \ref{series2}.

\vspace{.1in}
We finally present a direct proof of the fact that (\ref{ss}) is global in time for special
case $a_{02}=a_{03}=\cdots=0$.

\begin{thm}\label{sing} Consider (\ref{kdvp}) with $\nu$ and $\alpha$ satisfying $\nu^2+4\alpha^2\ge 9$. If
$$
u_0(x) =a_{01}\,e^{ix}\quad \mbox{with} \quad |a_{01}| < 1,
$$
then (\ref{kdvp}) has a unique global solution, which can be represented by (\ref{ss}).
In addition, for any $s\ge 0$, $u(\cdot,t)\in H^s$ for all $t\ge 0$.
\end{thm}
\noindent{\it Proof}. \,\, We prove by induction that, for any $t>0$,
\begin{equation}\label{ubb}
|a_k(t)| \le |a_{01}|^k, \quad k=1,2,\cdots.
\end{equation}
Obviously, $|a_1(t)|\le |a_{01}|$. To prove (\ref{ubb}) for $k\ge 2$, we recall (\ref{mit}), namely
$$
a_k(t) =  3 ik \,e^{-(\nu k^2 -\alpha i k^3)t}\int_0^t e^{(\nu k^2 -\alpha i k^3)\tau}
\sum_{k_1+k_2=k} a_{k_1}(\tau) \,a_{k_2}(\tau)\,d\tau.
$$
Since $\nu^2+4\alpha^2\ge 9$, we have
$$
|a_2(t)| \le \left|\frac{3}{2\nu -4\alpha i}\right| |a_{01}|^2\,
\left(1-e^{-(4\nu-8\alpha i)t}\right) \le |a_{01}|^2
$$
and more generally,
$$
|a_k(t)| \le \left|\frac{3(k-1)}{\nu\,k -\alpha i\,k^2}\right| |a_{01}|^k\,
\left(1-e^{-(\nu\,k^2-\alpha i\,k^3)t}\right) \le |a_{01}|^k.
$$
It is then clear that (\ref{ss}) converges in $H^s$ with $s\ge 0$ for any $t\ge 0$.
This completes the proof of Theorem \ref{sing}.

\vspace{.4in} \noindent {\bf Acknowledgements}

\vspace{.1in} NK and JW thank Professor David Wright for teaching
them the result of Hardy and Ramanujan and thank Professors Charles
Li and Ning Ju for discussions. JW also thanks Department of
Mathematical Sciences at University of Cincinnati for support and
hospitality. JMY thanks the NSC grant \# 96-2115-M126-001. BZ thanks
the Taft Memorial Fund  for its financial support.


\begin{thebibliography}{99}
\bibitem{Bir} B. Birnir,  An example of blow-up, for the complex KdV equation
and existence beyond the blow-up, {\it SIAM J. Appl. Math.  } {\bf 47} (1987), 710-725.

\bibitem{BW} J. L. Bona and F. B. Weissler, Blow up of spatially periodic complex-valued solutions
 of nonlinear dispersive equations, {\it Indiana Univ. Math. J.} {\bf 50} (2001), 759-782.

\bibitem{HR} G.H. Hardy and S. Ramanujan, Asymptotic Formulae in Combinatory
Analysis, {\it Proc. London Math. Soc. \bf 17} (1918), 75-115.

\bibitem{Ke}M. Kerszberg, A simple model for unidirectional crystal growth,
{\it Phys. Lett.} {\bf 105A} (1984), 241-244.

\bibitem{Le}D. Levi, Levi-Civita theory for irrotational water waves in a one-dimensional channel and the
complex Korteweg-de Vries equation, {\it Teoret. Mat. Fiz.} {\bf 99}
(1994), 435-440; {\it Translation in  Theoret. and Math. Phys.} {\bf
99} (1994), 705-709.

\bibitem{LeSa}D. Levi and M. Sanielevici, Irrotational water waves and
the complex Korteweg-de Vries equation, {\it Phys. D} \,{\bf 98} (1996), 510-514.

\bibitem{LS} D. Li and Y. Sinai, Blow ups of complex solutions of the 3D Navier-Stokes
equations, arXiv:Physics/0610101v1, 13 Oct 2006.

\bibitem{Li} Y. Charles Li, Simple explicit formulae for finite time blow up solutions to the complex KdV
equation, {\it Chaos, Solitons \& Fractals}, 2008 (in press).


\bibitem{PS} P. Pol\'{a}\v{c}ik and V. \v{S}ver\'{a}k, Zeros of complex caloric functions and singularities of
complex viscous Burgers equation, {\it J. Reine Angew. Math. \bf  616} (2008), 205--217.

\bibitem{Si}Y. Sinai, Power series for solutions of the $3D$-Navier-Stokes system on $R\sp 3$, {\it J. Stat. Phys.
  \bf 121}  (2005),  779--803.

\bibitem{WY1} J. Wu and J.-M. Yuan, The effect of dissipation on solutions of the complex KdV equation, {\it Math.
Comput. Simulation \, \bf 69}  (2005), 589--599.

\bibitem{WY2} J. Wu and J.-M. Yuan, Local well-posedness and local (in space) regularity results for the complex
Korteweg-de Vries equation, {\it Proc. Roy. Soc. Edinburgh Sect. A  \bf 137}  (2007), 203--223.

\bibitem{YW1} J.-M. Yuan and J. Wu, The complex KdV equation with or without dissipation,
{\it Discrete Contin. Dyn. Syst. Ser. B  \bf 5}  (2005), 489--512.



\end{thebibliography}
\end{document}